\input amstex
\documentstyle{amsppt}
\pagewidth{5.5in}

\magnification=1100


\def\ra{\rightarrow}
\def\lc{{{\Cal L}}}
\def\fg{{\frak g}}
\def\fw{{\frak w}}
\def\fv{{\frak v}}

\def\l2{L^2}
\def\rb{{\bold R}}

\def\Iso{\text{Iso}}

\def\trace{\text{trace}}

\def\Id{{\text{Id}}}
\def\ad{\text{ad}}
\def\so{\text{so}}

\def\span{\text{span}}

\def\Ric{\text{\rom{Ric}}}
\def\scal{\text{\rom{scal}}}
\def\proj{\text{proj}}

\def\jo{j_0}

\def\ogj{\overline{G(j)}}
\def\ogjp{\overline{G(j^\prime)}}

\def\fz{{\frak{z}}}

\def\bfz{\fz/\lc}

\def\nj{{N(j)}}
\def\njp{{N(j')}}

\def\tk{{T^k}}
\def\sn{{S^n}}

\def\ts{{\sn\times\tk}}
\def\nhj{{N_H(j)}}
\def\nhjp{{N_H(j')}}
\def\tp{{\tilde{p}}}
\def\tric{{\widetilde{\Ric}}}
\def\tscal{{\widetilde{\scal}}}
\def\bz{{\bar{z}}}
\def\nab{{\nabla}}
\def\tnab{{\widetilde{\nab}}}
\def\M{{\nj}}
\def\tM{{\ogj}}
\def\N{{\nu}}
\def\tang{{\text{T}_p\M}}
\def\ttang{{\text{T}_p\tM}}
 
\topmatter
 
\title Isospectral deformations of closed Riemannian manifolds
with different scalar curvature
\endtitle
 
\rightheadtext{Isospectral deformations with different scalar curvature}
\leftheadtext{Gordon, Gornet, Schueth, Webb, Wilson}
 
\author Carolyn S. Gordon, Ruth Gornet, Dorothee Schueth,
David L. Webb, Edward N. Wilson
\endauthor
 
\dedicatory
To the memory of Hubert Pesce \rom{(1966--1997)}, a gifted colleague and 
deeply missed friend.
\enddedicatory
 
\thanks
Much of this
work was carried out at a workshop in Grenoble in June 1997 funded by NSF-CNRS
grant INT-9415803 for international collaboration.  The authors also wish to
express their gratitude to
the Grenoble group for its fine organizational efforts.  The third author
wishes to
thank the GADGET program for supporting her participation.  
Research by the first and fourth
authors is partially supported by National Science Foundation grant
DMS-9704369. The second author is supported in part by the
Texas Advanced Research Program under Grant No. 003644-002.
The third author is supported in part by SFB 256, Bonn.
\endthanks



\abstract
We construct the first examples of continuous families of isospectral 
Riemannian metrics that are not locally isometric on closed manifolds , more
precisely, 
on $S^n\times T^m$, where $T^m$ is a torus of dimension $m\ge 2$ and $S^n$ is a
sphere 
of dimension $n\ge 4$.  These metrics are not locally homogeneous;  in
particular, the 
scalar curvature of each metric is nonconstant.  For some of the deformations,
the 
maximum scalar curvature changes during the deformation.
\endabstract

\address
\newline Carolyn S. Gordon, Dartmouth College, Hanover, NH 03755-3551,
\newline \indent carolyn.s.gordon\@dartmouth.edu
\newline Ruth Gornet, Texas Tech University, Lubbock, TX 79409-1042,
\newline \indent gornet\@math.ttu.edu
\newline Dorothee Schueth, Mathematisches Institut, Universit\"at Bonn,
Beringstr\. 1, D-53115 Bonn,
\newline\indent schueth\@math.uni-bonn.de
\newline David L. Webb, Dartmouth College, Hanover, NH 03755-3551,
\newline \indent david.l.webb\@dartmouth.edu
\newline Edward N. Wilson, Washington University, St\. Louis, MO 63130-4899,
\newline \indent enwilson\@math.wustl.edu
\endaddress
 
\endtopmatter

\document
 
\subheading{Introduction}
\smallskip
 
To what extent does the eigenvalue spectrum of a compact Riemannian manifold
determine the geometry of the manifold?  Various global geometric invariants,
such as dimension, volume, and total scalar curvature, are known to be
spectrally determined.
Moreover, various manifolds such as round spheres of dimension less than or
equal to six and $2$-dimensional flat tori are uniquely determined by their
spectra. However, we show:
 
\proclaim{Main Theorem} For $n\geq 4$, there exist
continuous $d$-parameter families $\{g_t\}$ of isospectral, non-isometric
Riemannian metrics on the manifold $S^n \times T$, where
$T$ is the $2$-dimensional torus and $\sn$ is the $n$-dimensional sphere.
Here $d$ is of order at least $O(n^2)$. These metrics are non-homogeneous.
For some of the deformations, the maximum scalar curvature of $g_t$ depends
non-trivially on $t$.\endproclaim
 
The metrics can be chosen arbitrarily close to the standard metric, i.e., to
the product of the flat metric on the torus and the round metric on the
sphere.
 
The first example of isospectral Riemannian manifolds was a pair of sixteen
dimensional flat tori \cite{M}. The past dozen years have seen an explosion
of new examples. See, for example, \cite{BT}, \cite{Bu},
\cite{GWW}, \cite{GW1,2}, \cite{Gt1,2}, \cite{I}, \cite{Sch}, \cite{Su},
\cite{V}.
However, until recently, all known isospectral manifolds
were at least locally isometric; in particular, all isospectral closed
manifolds had a common Riemannian cover. This was due primarily to the
fact that most examples could be explained by Sunada's method \cite{Su} or
its generalizations \cite{DG2}. The Sunada methods rely almost exclusively
on representation theory, with the result that isospectral manifolds
constructed using these methods must be locally isometric. See the
expository articles \cite{Be}, \cite{Br}, \cite{D}, \cite{G1}, or \cite{GGt}
for more
information about isospectral manifolds in general.
 
Then Szabo \cite{Sz} constructed pairs of isospectral compact manifolds with
boundary that are not locally isometric, and Gordon and Wilson \cite{GW3}
generalized his construction to obtain
continuous families of such manifolds. Finite families of closed
isospectral manifolds with different local geometry were given in
\cite{G2,3} and \cite{GW3}.
 
The examples described in the Main Theorem above have several new features:
\roster
\item"{$\bullet$}" They give the first examples of continuous isospectral
deformations of closed manifolds for which the metrics are not locally
isometric.
\item"{$\bullet$}" They are essentially the first examples of isospectral
metrics that are not locally homogeneous, e.g., for which the curvature
varies from point to point. (The only other examples are locally isometric
metrics
on nilpotent Lie groups that are not fully left-invariant \cite{DG1} and
product (or twisted products)
$M\times N$ and $M'\times N$ where $M$ and $M'$ are locally homogeneous
isospectral manifolds and $N$ is fixed.)
\item"{$\bullet$}" They are the first examples of isospectral manifolds with
different scalar curvature.\endroster
 
The isospectral manifolds in the Main Theorem are the boundaries of the
isospectral manifolds constructed in \cite{GW3}.
\smallskip

After completing a draft of this paper, the authors learned that Z. Szabo has 
independently shown that the boundaries of the manifolds he constructed in
\cite{Sz} 
are isospectral.  Like the isospectral deformations considered here, the pairs
of isospectral 
metrics that he constructed are not locally homogeneous nor locally isometric,
although they do have the same 
maximum and minimum scalar curvature.  Szabo, moreover, explicitly computed the
spectrum 
of these isospectral metrics.  This work is included in a revised version of
the article 
\cite{Sz}.

\subheading{Construction of isospectral metrics}
\smallskip
 
The isospectral manifolds constructed in \cite{GW3} are domains in certain
two-step nilpotent Lie groups. They may be described as follows.
 
\definition{Notation 1} Let
$\fz=\rb^k$ and $\fv=\rb^m$ with their standard inner products, and let
$j:\fz\to\so(\fv)$ be a linear map. Let $\fg$ be the orthogonal
direct sum $\fg =\fv\oplus\fz$. Define a Lie bracket on $\fg$ by declaring
$\fz$ to be central and defining $[\,,\,]:\fv\times\fv\to\fz$ by

$$\left<[x,y],z\right>=\left<j(z)x,y\right>
$$
for all $x$ and $y$ in $\fv$ and $z$ in $\fz$. This bracket gives $\fg$ the
structure of a two-step nilpotent Lie algebra, and
all two-step nilpotent Lie algebras may be constructed by this method.
We denote $\fg$ with this structure by $\fg(j)$. An element $x\in\fv$ is
central if and only if $j(z)x=0$ for all $z\in\fz$.
Without loss of generality, we assume $\fv$ intersects the center trivially;
i.e., $\fz$ coincides with the center of $\fg(j)$.
(Note that $\fg(j)$ may still have an abelian factor,
as we are not assuming that $j$ is non-singular.)
Let $G(j)$ denote the simply-connected Lie group with Lie algebra
$\fg(j)$. The Lie group
exponential map $\exp:\fg(j)\to G(j)$ is a diffeomorphism. Moreover, the
center of $G(j)$ is isomorphic to $\rb^k$ and $\exp_{|\fz}$ is a linear
isomorphism from $\fz$ to the center of $G(j)$.
 
With respect to the global coordinate system
$G(j)\to\fv\oplus\fz$ defined by $\exp(x+z)\mapsto(x,z)$ for 
$x\in\fv$ and $z\in\fz$, the Lie group multiplication is given by

$$(x,z)\cdot (x',z')=(x+x',z+z'+\frac{1}{2}[x,y]).
$$

The inner product on $\fg(j)$ defines a left-invariant metric on $G(j)$,
i.e., a metric for which the left translations by group elements are
isometries.
 
Let $\lc$ be a lattice of
full rank in $\fz$. Identify $\lc$ with $\exp(\lc)$, a discrete central
subgroup of $G(j)$.
The quotient $\ogj=G(j)/\lc$ is again a two-step nilpotent Lie group,
and it inherits a Riemannian metric from its covering $G(j)$. Its center
is isomorphic to the $k$-dimensional torus $\fz/\lc$, and
$\ogj$ is diffeomorphic to $\fv\times (\fz/\lc)$. Let $B$ be the closed unit
ball in $\fv$ and $S$ its $(m-1)$-dimensional boundary sphere. Set

$$M(j)=\{(x,\bz)\in\ogj : x\in B\text{ and }\bz\in\fz/\lc\}
$$

and

$$N(j)=\{(x,\bz)\in\ogj : x\in S\text{ and }\bz\in\fz/\lc\}.
$$

\enddefinition
 
The manifolds $M(j)$ were studied extensively in \cite{GW3}.
The $M(j)$ are manifolds with boundary that are locally homogeneous;
their local geometry is that of the Lie group $G(j)$.
We study their boundaries $N(j)$ here. Note that
$N(j)$ is diffeomorphic to $\ts$ where $T^k$ is a torus of dimension
$k=\dim(\fz)$ 
and where $n=m-1$.
 
\proclaim{Proposition 2} The map $\pi_j :N(j)\to S$ given by
$(x,\bz)\to x$ is a Riemannian submersion with respect to the
canonical round metric on $S$ associated with the inner product on $\fv$.
The fibers are totally geodesic flat tori that are orbits of a flat toral
group $T(j)\simeq \fz/\lc$ of isometries acting freely on $N(j)$.
\endproclaim
 
\demo{Proof}
The analogous statement for the nilmanifold $\ogj$ is well-known; i.e., the
map $\ogj\to \fv$ given by $(x,\bz)\mapsto x$ is a Riemannian submersion with
respect to the Euclidean metric on $\fv$, and the fibers are totally geodesic
flat tori. The fibers are orbits of the toral group of isometries given by
translations by the elements $(0,\bz)$, $\bz\in\bfz$. The manifold
$N(j) \subset \ogj$ is the restriction of the torus bundle $\ogj$ over $B$ to
the
subspace $S \subset B$, and the statement follows.
\qed\enddemo
 
\definition{Definition 3}  Let $\fv$, $\fz$ and $\lc$ be as above.
 
\item{(i)} A pair $j$, $j'$ of linear maps from
${\fz}$ to $\so({\fv})$ is called
$\lc$-{\it equivalent} if there exist orthogonal linear operators $A$ on
$\fv$  and $C$ on $\fz$ such that $C(\lc)=\lc$ and

$$Aj(z)A^{-1} = j'(C(z))
$$

for all $z \in \fz$. We call the pair $(A,C)$ an
$\lc$-{\it equivalence} between $j$ and $j'$.
 
\item{(ii)} We say $j$ is {\it isospectral} to $j'$, denoted
$j \sim j'$, if for each $z \in \fz$, the eigenvalue spectra (with
multiplicities) of $j(z)$ and $j'(z)$ coincide, i.e., there exists an
orthogonal linear operator $A_z$ for which

$$ {A_z j(z) A_z^{-1} = j'(z)}.
$$

\enddefinition
 
\subheading{Remark 4} We address the question of when $N(j)$ is
isometric to $\njp$ later in the paper. For now, we note the following
facts about $\ogj$ and $\ogjp$.
\item{(i)} Given orthogonal transformations $A$ of $\fv$ and $C$ of $\fz$
such that $C(\lc)=\lc$,
one easily checks that the following statements are equivalent:

\roster

\item"{(a)}" The pair $(A,C)$ is an $\lc$-equivalence between $j$ and $j'$.

\item"{(b)}" The map
$\fg(j)\to\fg(j')$ given by $x+z\mapsto A(x)+C(z)$ for $x\in\fv$ and
$z\in\fz$ is a Lie algebra isomorphism.

\item"{(c)}" Viewing $C$ as an automorphism of the torus $\bfz$, the map
$\tau:\ogj\to\ogjp$ given by $\tau(x,\bz)= (A(x),C(\bz))$ is both a Lie group
isomorphism and an isometry.

\endroster
 
\item{(ii)} By \cite{W},
any isometry between $\ogj$ and $\ogjp$ is of the form
$L_{(x,\bz)}\circ\tau$ where $\tau:\ogj\to\ogjp$ is both a Lie group
isomorphism and an isometry (and thus arises from an $\lc$-equivalence
between $j$ and $j'$)
and $L_{(x,\bz)}:\ogjp\to\ogjp$ is left translation by the element
$(x,\bz)$ of $\ogjp$. In particular, $\ogj$ is isometric to $\ogjp$ if and
only if $j$ and $j'$ are $\lc$-equivalent.
 
\item{(iii)} If $j$ is $\lc$-equivalent to $j'$,
then the isometry $\tau:\ogj \to\ogjp$ defined in (i)(c)
restricts to an isometry from $N(j)$ to $N(j')$.

\proclaim{Theorem 5}  Fix inner product spaces $\fv$ and $\fz$ and a lattice
$\lc$ of full rank in $\fz$, and let $j,j':\fz\to\so(\fv)$ be linear maps.
If $j\sim j'$, then $N(j)$ is isospectral to $N(j')$.
\endproclaim
 
\demo{Proof} The analogous statement for $M(j)$ and $M(j')$
is proved in \cite{GW3}. The proof that
$N(j)$ and $N(j')$ are isospectral is similar, so we just
give the main ideas here. First note that if
$\fz$ is one-dimensional and if $j\sim j'$, then $j$ is $\lc$-equivalent to
$j'$, with $C$ in Definition 3  being the identity map.
By Remark 4, the manifolds are isometric, hence trivially isospectral.
 
Now consider the general case. By Proposition 2, the central torus
$T(j)\simeq\fz/\lc$ acts on $N(j)$ by isometries, and the manifold
$N(j)$ has the structure of a torus bundle over the round sphere $S$.
Moreover, the projection $N(j)\to S$ is a
Riemannian submersion with totally geodesic fibers. If $\fw$ is a
codimension-one subspace of $\fz$ spanned by lattice vectors in $\lc$, then
$\fw$ covers a codimension-one subtorus $H$ of $T(j)$. Let $N_H(j)$ denote the
quotient of
$N(j)$ by the isometry action of this torus. With the induced metric on
$\nhj$, the projection
$N(j)\to\nhj$ is again a Riemannian submersion with totally geodesic fibers.
 
Now $\nhj$ is defined in the same way as $N(j)$, but with $\fz$ replaced by
$\fw^\perp =\fz\ominus\fw$, the orthogonal complement of
$\fw$ in $\fz$, and $j$ replaced by $j_{|{\frak w}^\perp}$. Since
$j_{|{\frak w}^\perp}\sim j'_{|{\frak w}^\perp}$ and since ${\frak w}^\perp$ is 
one-dimensional, the first part of the proof shows that $\nhj$ is isometric to
$\nhjp$.
 
Functions on $\nhj$ pull back to functions on $N(j)$. Since the projection is a
Riemannian submersion with totally geodesic fibers, a classical theorem says
that eigenfunctions on $\nhj$ pull back to eigenfunctions on $N(j)$ with the
same
eigenvalue. Thus the spectrum of $\nhj$ is part of the spectrum of $N(j)$.
By a Fourier decomposition argument involving the action of the torus $T(j)$,
one shows that the spectra of the various $\nhj$
exhaust the spectrum of $N(j)$.
We can then use the isospectrality of the $\nhj$ with the $\nhjp$ to
deduce the isospectrality of $N(j)$ with $N(j')$.
\qed\enddemo
\smallskip
 
\subheading{Scalar Curvature}
\smallskip
 
We next express the scalar curvature function $\scal$ of the manifold $N(j)$
in terms of the scalar and Ricci curvatures $\tscal$ and $\tric$ of the
ambient nilmanifold $\ogj$. Since $\ogj$ is locally homogeneous,
$\tscal$ is a constant function; we denote its constant value also by
$\tscal$. The Ricci curvature $\tric$ does not depend on the point and may
be viewed as a bilinear form on the Lie algebra $\fg(j)$.
 
\proclaim{Proposition 6} Using Notation 1, let $m=\dim(\fv)$.
The scalar curvature $\scal$ of
$N(j)$ is given at the point $p=(x,\bz)$ of $N(j)$ by

$$\scal(p)=\tscal+(m-1)(m-2)-\tric(x,x).
$$

\endproclaim
 
\demo{Proof}
Let $\N$ denote a unit normal vector field to $\M$ in $\tM$.
Denote by $\tnab$, $\nab$ the Levi--Civita connections of
$\tM$ and $\M$ respectively.
 
A routine calculation using the second fundamental
form shows that for $p$ in $\M$,

$$\scal(p)=\tscal-2 \tric(\N_p,\N_p)+(\trace(\tnab \N_{\,|\,\tang} ))^2
-\Vert\tnab \N_{\,|\,\tang}\Vert^2.\tag6.1
$$
(The norm in the last term is that on tensors of type (1,1), i.e., 
the $L^2$ norm on matrices relative to an orthonormal basis.)

We remark that for any two-step nilpotent Lie group $G$ with Lie algebra $\fg$ 
and any $v$ in $\fg$, one has the identity

$$\exp_{*v}=L_{\exp(v)*}(\Id-\frac1{2}\ad_v).                           
\tag6.2
$$

Let $p=(x,\bz)$ and recall that $\ttang=L_{p*}(\fg(j))$, where $L_p$ is
left translation by $p$.
Using (6.2), one easily checks that as a subspace of $L_{p*}(\fg(j))$,
$\tang=L_{p*}(x^\perp)$,
where $x^\perp$ is the orthogonal complement of $\span\{x\}$ in $\fg(j)$.
The outward unit normal vector field is then given by  $\nu_p=L_{p*}(x)$.
 
Let $y\in x^\perp \cap \fv$ and $w\in\fz$.
We view $x$, $y$, $\tnab_y x$, $w$, and $\tnab_w x$ as left-invariant vector
fields on $\ogj$, so we write, e.g., $y_p$ for $L_{p*}y$.
The curve $c(t)=\exp((\cos t)x+(\sin t)(y+\frac1{2}[x,y])+\bz)$ has initial
velocity $y_p$.
Viewing $\N$ as a vector field along $c$ and computing the covariant derivative
yields

$$\tnab_{y_p} \N=y_p+(\tnab_y x)_p.                                     
$$
Similarly, 
$$     
\tnab_{w_p} \N=(\tnab_w x)_p.        										$$

In summary, for $p=(x,\bz) \in\M$ and $u\in x^\perp=L_{p*}^{-1}(\tang)$,

$$
{\tnab}_u \N =L_{p*}(\proj_\fv u+\tnab_u x),
$$
where $x$ is viewed as a left-invariant vector field.
 
Recall (see \cite{E, Section 2}) that for $x,y\in\fv$ and $w\in\fz$,

$$\tnab_y x = \frac1{2}[y,x]\quad\text{and}\quad
\tnab_w x =-\frac1{2}j(w)x.\tag6.3
$$

Thus $\tnab x$ sends $\fv\to\fz$ and $\fz\to\fv$, and
it follows immediately that $\trace(\tnab x)=0$, and that
$\trace(\tnab\N_{\,|\,\tang})=m-1$.
Likewise, $\Vert(\tnab \N)_{\,|\,\tang}\Vert^2=(m-1)+\Vert\tnab x\Vert^2$.
 
Finally, let $z_1,\dots,z_k$ be an orthonormal basis of $\fz$.
It is known that (see \cite{E, Section 2})

$$\tric(x,x)=\frac1{2}\sum_{i=1}^k\left<j(z_i)^2 x , x \right>.\tag6.4
$$

Equations (6.3) and (6.4) and Notation 1 imply

$$\Vert\tnab x\Vert^2=-\tric(x,x)
$$
and Proposition 6 follows.
\qed\enddemo
 
\proclaim{Corollary 7}  Fix inner product spaces $\fv$ and $\fz$ and
a lattice $\lc$  of full rank in $\fz$,  and let $j,j':\fz\to\so(\fv)$
be linear maps. The maximum, respectively minimum, values of the scalar
curvature of $N(j)$ and $N(j')$ coincide if and only if the minimum,
respectively maximum, eigenvalues of the Ricci tensors of $\ogj$ and
$\ogjp$ restricted to $\fv\otimes\fv$ coincide.
\endproclaim
 
\subheading{Examples 8} We briefly review Example 2.3 of \cite{GW3}
and show that some of the isospectral deformations constructed there have
changing maximal scalar curvature.
 
Take $\fz=\rb^2$ and $ \fv=\rb^6$ with their standard ordered bases
and standard inner product. For $a,b\in\so(6)$ and $s,t\in\rb$, define
$j_{a,b}(s,t)=sa+tb$. Each linear map
$j:\rb^2\ra\so(6)$ is of the form $j=j_{a,b}$ for some $a,b\in \so(6)$.
Fix for the remainder of the discussion an element $a\in\so(6)$ that is in
block diagonal form with $2\times 2$ diagonal blocks
$\left[\matrix 0&-a_i\\
a_i&0\endmatrix\right], 1 \leq i
\leq 3$, where $0<a_1<a_2<a_3$.
Consider all matrices $b\in \so(6)$ of the form

$$ b=\left[\matrix
0&0&b_{12}&0&b_{13}&0\\
0&0&0&0&0&0\\
-b_{12}&0&0&0&b_{23}&0\\
0&0&0&0&0&0\\
-b_{13}&0&-b_{23}&0&0&0\\
0&0&0&0&0&0\endmatrix\right]
$$
with $(b_{12},b_{13},b_{23}) \in\rb^3$.
From \cite{GW3}, we know that if $b$ and $b'$ are of this form, then
$j_{a,b}\sim j_{a,b'}$ if and only if there exists a real number
$u$ in $I=[\max\left(\frac{-b_{12}^2}{a_2^2-a_1^2},\frac{-b_{23}^2}
{a_3^2-a_2^2}\right),\frac{b_{13}^2}{a_3^2-a_1^{2}}]$
satisfying

$$\align (b_{12}^\prime)^2&=b_{12}^2 +u(a_2^2-a_1^2),\\
(b_{13}^\prime)^2&=b_{13}^2
+u(a_1^2-a_3^2),\tag$*$\\ (b_{23}^\prime)^2&=b_{23}^2+u(a_3^2-a_2^2).
\endalign
$$

If we fix any $b$ for which $I$ has non-empty interior and,
for each $u\in I$, define a one-parameter family 
$b(u)$ as the unique solution $b^\prime$ of the above equations for which
$b_{ij}(u)$ has the same sign as $b_{ij}$ for all $i,j$, it follows
that $u\ra j_{a,b(u)}$ is a $1$-parameter isospectral deformation of
$j_{a,b}$.
 
Using (6.4) above, we see that
$\tric(u)(x,x)=\frac1{2}\left<(a^2+b(u)^2)x, x\right>$ for $x\in\fv$.
For example, letting $a_1=1$, $a_2=2$, $a_3=3$, $b_{12}=0$, $b_{13}=1$,
and $b_{23}=0$, we obtain, with respect to the standard basis on $\fv$,

$$\tric(u)= -\frac1{2}\left[\matrix
2-5u&0&\sqrt{5u-40u^2}&0&-\sqrt{15}u&0\\
0&1&0&0&0&0\\
\sqrt{5u-40u^2}&0&4+8u&0&\sqrt{3u-24u^2}&0\\
0&0&0&4&0&0\\
-\sqrt{15}u&0&\sqrt{3u-24u^2}&0&10-3u&0\\
0&0&0&0&0&9\endmatrix\right]
$$
for $u\in[0,\frac1{8}]$.
One easily checks that the eigenvalues of $\tric(u)$ change
non-trivially with $u$; in particular, $-5$ is an eigenvalue
of $\tric(0)$, but the eigenvalues of $\tric(u)$ have absolute 
values strictly less
than $5$ for $0<u<\frac1{8}$. By Corollary 7, the maximum value of the
scalar curvature of $N(j_{a,b(u)})$ changes non-trivially
with $u$, showing that $N(j_{a,b(u)})$ is a non-trivial isospectral
deformation.
\smallskip
 
\subheading{Large-Dimensional Families}
\smallskip
 
Finally we turn to the Main Theorem, stated in the introductory remarks. In
the proposition below, the notion of {\it equivalence} of linear maps
$j,j':\fz\to\fv$ is identical to the notion of $\lc$-equivalence in
Definition 3 except that the map $C$ is not required to preserve a given
lattice $\lc$.
 
\proclaim{Proposition 9 \cite{GW3}}
Let $\dim\fz=2$, and let $m=\dim {\frak v}$ be any positive integer other
than $1,2,3,4$, or $6$. Let $W$ be the real vector space consisting of
all linear maps from $\fz$ to $\so({\frak v})$. Then there is a Zariski open
subset ${\Cal O}$ of $W$ (i.e., ${\Cal O}$ is the complement of the zero locus 
of some non-zero polynomial
function on $W$) such that each
$\jo\in {\Cal O}$ belongs to a $d$-parameter family of isospectral,
inequivalent elements of $W$. Here $d\geq m(m-1)/2 - [m/2]([m/2]+2)>1$.
In particular, $d$ is of order at least $O(m^2)$.
\endproclaim
 
Although the expression for $d$ gives $0$ when $m=6$, Example 8 gives
continuous families of isospectral, inequivalent $j$ maps when $m=6$.
 
Letting $n=m-1$, choosing any lattice $\lc$ in $\fz$, and applying Theorem 5,
the families of isospectral $j$ maps in Proposition 9 give rise to continuous
families of isospectral manifolds $N(j)$, or equivalently
to continuous families of isospectral metrics on
$S^n\times T$ where $T\simeq \fz/\lc$ is a $2$-torus and $n=m-1$. To
complete the proof of the Main Theorem, we need to show that the manifolds in
these deformations are not isometric. Recall that by Remark 4,
$\ogj$ and $\ogjp$ are isometric if and only if $j$ and $j'$ are
$\lc$-equivalent. We now address the question of whether $\lc$-equivalence
of $j$ and $j'$ is necessary in order for the submanifolds $N(j)$ and
$N(j')$ to be isometric.
 
\proclaim{Proposition 10}
Fix inner product spaces $\fv$ and $\fz$ and a lattice
$\lc$ of full rank in $\fz$,  and let $j,j':\fz\to\so(\fv)$ be linear maps.
Suppose that there are only finitely many automorphisms of $\fg(j)$ that
restrict to the identity on $\fz$. If $N(j)$ is isometric to $N(j')$, then
$j$ is $\lc$-equivalent to $j'$.
\endproclaim

\subheading{Remark 11} For any choice of $j$,
the linear map that restricts to $-\Id$ on $\fv$ and to $\Id$ on
$\fz$ is an automorphism of $\fg(j)$. For generic choices of $j$, this is the
only non-trivial automorphism that restricts to the identity on $\fz$.
\smallskip 

Before proving Proposition 10, we note that the Main Theorem now
follows immediately:

\demo{Proof of Main Theorem}
While the proof of Proposition 9 in \cite{GW3} is not constructive,
the only maps $j\in W$ considered there are those that satisfy the genericity
condition in Remark 11, hence the hypothesis of Proposition 10.
Thus each $d$-parameter family $\{j_t\}$ of isospectral, inequivalent
$j$ maps gives rise to a $d$-parameter family of isospectral, non-isometric
manifolds $N(j_t)$.
\qed\enddemo
 
In order to prove Proposition 10, we need the following Lemma.

\proclaim{Lemma 12}  Fix inner product spaces $\fv$ and $\fz$ and a lattice
$\lc$  of full rank in $\fz$,  and let $j,j':\fz\to\so(\fv)$ be linear maps.
If $\tau:N(j)\to N(j')$ is both an isometry and a
bundle map with respect to the bundle structures defined in Proposition 2,
then $\tau$ extends to a map $\tilde{\tau}:\ogj\to\ogjp$ that is both a Lie
group isomorphism and an isometry.
Moreover, $j$ is $\lc$-equivalent to $j'$.
\endproclaim

\demo{Proof of Lemma 12}  Since $\tau$  preserves the fiber
structures, it induces an isometry on the unit sphere $S$ in $\fv$.
Such an isometry is given by an orthogonal transformation $A$ of $\fv$.
 
Since the principal torus bundles $N(j)$ and $N(j')$ over $S$ are topologically 
trivial, we can choose global sections of $N(j)$ and $N(j')$. The
torus action then determines an isomorphism of each fiber of
each bundle with the torus
$T=\fz/\lc$. For each
$p\in S$, the isometry $\tau$ must restrict to an isometry $\tau_p$ of the
fiber $T_p$ above $p$ in $N(j)$ to the fiber $T_{A(p)}$ above $A(p)$ in
$N(j')$. Now
$\tau_p$ must be of the form $C_p\circ L_{\bz_p}$, where
$C_p$ is an orthogonal transformation of $\fz$ that preserves $\lc$, and
$L_{\bz_p}$ denotes translation by some element $\bz_p$ of the torus $T$.
Both $\bz_p$ and $C_p$ must vary continuously with $p$.
However, $C_p$ stays within a discrete set and hence must be independent
of $p$; we thus drop the subscript.
 
We now show that $\ogj$ and $\ogjp$ must be isomorphic. Let $x$ and $y$ be
any orthonormal pair of vectors in $\fv$. The $2$-plane spanned by these
vectors intersects $S$ in a great circle. Choose a lift $\tp$ of $x$ in $N(j)$. 
Lift the great circle to a
horizontal geodesic $\sigma$ in $N(j)$ starting at $\tp$. Since $\tau$
carries the vertical space at $\tp$ to the vertical space at $\tau(\tp)$,
it must induce an isomorphism between the horizontal spaces as well; this
isomorphism is determined by the orthogonal transformation $A$ of $\fv$.
In particular,
$\tau$ must carry $\sigma$ to a horizontal geodesic  through $\tau(\tp)$;
this horizontal geodesic is the lift of the great circle defined by the
$2$-plane of the vectors $A(x)$ and $A(y)$.
 
If $\tp=(x,\bz)$,
the geodesic $\sigma$ is given by $\sigma(t)=((\cos t)x+(\sin t)y, \bz(t))$
with $z(t)=z+t[x,y]/2$. Thus $\sigma$ in general is not closed
unless $x$ and $y$ commute. The displacement in the torus fiber between the
initial and final points is given by

$$
\sigma(2\pi)-\sigma(0)=\pi[x,y].
$$

Similarly,

$$(\tau\circ\sigma)(2\pi)-(\tau\circ\sigma)(0)=\pi[A(x),A(y)].
$$

On the other hand,

$$(\tau\circ\sigma)(2\pi)-(\tau\circ\sigma)(0)=C(\sigma(2\pi)-\sigma(0)).
$$
These equations are to be understood in the sense that all terms are
well-defined modulo~$\lc$. 

We conclude that

$$C([x,y])\equiv[A(x),A(y)] \ (\text{mod}\ \lc) \tag12.1$$
for each pair of orthonormal vectors $x,y\in\fv$. We claim that the equivalence
is 
actually an equality, in fact, that

$$C([x,y])=[A(x),A(y)] \tag12.2
$$
for all $x,y\in\fv$.  If $\dim(\fv)\ge 3$, the set $\Cal P$ of pairs
of 
orthonormal vectors forms a connected subset of $\fv \times \fv$.  By
continuity
of $A$ and $C$, there must exist a constant vector $a\in\lc$ such that 
$C([x,y])-[A(x),A(y)]=a$ for all $(x,y)\in {\Cal P}$. Replacing $x$ by
$-x$ in this equation, we see that $a$ must be zero, i.e., that
$$C([x,y])=[A(x),A(y)]$$
for all $(x,y)\in {\Cal P}$ and thus for arbitrary vectors $x,y\in\fv$, as
claimed.  
While the case $\dim(\fv) < 3$  is not of interest in our application of
Lemma 12,
we note that equation (12.2) can also be proven in all dimensions by lifting
the 
isometry $\tau$ to an isometry between the universal coverings of $\nj$ and
$\njp$ and carrying out the
entire proof leading up to the equivalence (12.1) at this level.  

From equation (12.2), it follows that the linear map from
$\fg(j)$ to
$\fg(j')$ given by $x+z\mapsto A(x)+C(z)$ for $x\in\fv$ and $z\in\fz$ is a
Lie algebra isomorphism as well as an orthogonal map. By Remark~4, the pair
$(A,C)$ is an $\lc$-equivalence between $j$ and $j'$, and thus
$\ogj$ is isometric to $\ogjp$.

However, we still need to show that $\tau$ itself extends to an isometry
$\tilde{\tau}:\ogj\to\ogjp$. Let
$\alpha:\ogj\to\ogjp$ be the isometry associated with the $\lc$-equivalence
as in Remark 4(i). Then
$\alpha_{|N(j)}:\nj\to\njp$ is an isometry. Moreover,
$\zeta:=\alpha_{|\nj}^{-1}\circ\tau$ is an isometry of
$N(j)$ that preserves every fiber and acts on the fiber over $p\in S$ as
translation by $\bz_p$.
 
Fix a point $q\in S$. Compose $\zeta$ with an
appropriate element of $T(j)$ to obtain a fiber-preserving isometry $\beta$
of $N(j)$ that acts trivially on the fiber over $q$. Fix a point $u$ in this
fiber. Then the differential
$\beta_{*u}$ restricts to the identity on the tangent space to the fiber.
Moreover, since $\beta$ induces the identity on $S$, $\beta_{*u}$ also
restricts to the identity on the horizontal space at $u$. Hence we have both
$\beta(u)=u$ and $\beta_{*u}=\Id$. Since any isometry is uniquely determined
by its value and its differential at a single point, it follows that
$\beta =\Id$ and so $\zeta\in T(j)$. Thus $\zeta$ extends to a left translation
of $\ogj$ and $\tau=\alpha\circ\zeta$ extends to an isometry of $\ogj$.
\qed\enddemo
 
\demo{Proof of Proposition 10} Let $\Iso(j)$ be the identity component in the
full isometry group of $N(j)$. Since $N(j)$ is compact,
$\Iso(j)$ is a compact Lie group containing the toral group $T(j)$. We first
show that $T(j)$ is a maximal torus in $\Iso(j)$. Any isometry
$\alpha$ that commutes with
$T(j)$ must preserve the fiber structure. By Lemma 12, $\alpha$ extends to
an isometry $\tilde{\alpha}:\ogj\to\ogj$. By Remark 4(ii),
$\tilde{\alpha}=L_{(0,\bz)}\circ\beta$ for some $\bz\in\bfz$ and some
automorphism $\beta$ of $\ogj$. The left translation $L_{(0,\bz)}$ is the
extension to $\ogj$ of an element of $T(j)$. Moreover, since
$\alpha$ commutes with
$T(j)$, the automorphism $\beta$ commutes with translation by all central
elements of $\ogj$, and thus $\beta$ restricts to the identity on the center
of $\ogj$. By the hypothesis of the proposition, $\beta$ lies in a finite
set of automorphisms. Hence $T(j)$ must have finite index in its
centralizer, so $T(j)$ is a maximal torus in $\Iso(j)$.
 
Now suppose $\tau :N(j)\to N(j')$ is an isometry. Then the map
$\hat{\tau}:\Iso(j)\to \Iso(j')$ given by $\hat{\tau}(\beta)=
\tau\beta\tau^{-1}$ is an isomorphism from $\Iso(j)$ to $\Iso(j')$ and thus
carries $T(j)$ to a maximal torus in $\Iso(j')$. Since
$\dim(T(j))=\dim(\fz)=\dim(T(j'))$, it follows that $T(j')$ is a maximal
torus in $\Iso(j')$. Since all maximal tori in a compact Lie group are
conjugate, we may assume, after composing $\tau$ with an isometry of $N(j')$,
that $\hat{\tau}(T(j))=T(j')$. It follows that $\tau$ is a bundle map, and
by Lemma 12, $j$ is $\lc$-equivalent to $j'$.
\qed\enddemo
 
\subheading{Remark 13} Under the hypothesis of the Proposition,
we actually have that
$\Iso(j) = T(j)$. To see this, fix a point $u\in\nj$, and define
$\rho:\Iso(j)\to \nj$ by $\rho(\alpha)=\alpha(u)$. Then the restriction
$\rho_{|T(j)}$ maps $T(j)$ to a fiber in $\nj$ and induces an isomorphism on
fundamental groups
$\rho_*:\pi_1(T(j))\tilde\to\pi_1(N(j))$. On the other hand, the
inclusion $\iota:T(j)\to \Iso(j)$ of the maximal torus into the compact
Lie group $\Iso(j)$ induces an injection on fundamental groups only if
$\Iso(j)=T(j)$. (Indeed, a compact Lie group $K$ is a product of a compact
semisimple subgroup $H$ and a torus. Any maximal
torus $T$ must include a maximal torus in $H$. Since $H$ has finite
fundamental group, the inclusion of a maximal torus in $K$ cannot induce an
injection of fundamental groups unless $H$ is trivial. On the other hand, if
$H$ is trivial, then $K$ is itself a torus and thus must equal $T$.)  Since
$\rho\circ\iota =\rho_{|T(j)}$, we conclude that $T(j)=\Iso(j)$.

\enddocument